\documentstyle{article}

\input{amssym.def}
\input{amssym.tex}
\begin{document}
\title{Thurston's norm revisited}

\author{Igor V. Nikolaev
\footnote{Partially supported 
by NSERC.}}

\date{}
\maketitle


\newtheorem{thm}{Theorem}
\newtheorem{lem}{Lemma}
\newtheorem{dfn}{Definition}
\newtheorem{rmk}{Remark}
\newtheorem{cor}{Corollary}
\newtheorem{prp}{Proposition}
\newtheorem{con}{Conjecture}

\begin{abstract}
We study the Thurston  norm on the second homology
of a 3-manifold $M$, which is the surface bundle over the circle with a pseudo-Anosov 
monodromy.  A novelty  of our approach consists in the application of the $C^*$-algebras to
a problem in topology.  Namely, one associates to $M$  a $C^*$-algebra, 
whose $K$-theory gives rise to an algebraic number
field $K$. It is shown,  that the trace function on the ring
of integers of $K$ induces a norm on the second homology
of $M$. The norm coincides  with the Thurston norm on  
the second homology of $M$.

\vspace{7mm}


\vspace{5mm}
{\it AMS (MOS) Subj. Class.:  19K, 46L, 57M.}
\end{abstract}

\section{Introduction}
In 1986,  W.~P.~Thurston discovered the fundamental measure
of complexity of a $3$-dimensional manifold, $M$. 
The measure  is a norm on  the second homology of $M$, which assigns the non-negative integers to 
the elements of the group  $H_2(M)=H_2(M;{\Bbb Z})$. If $z\in H_2(M)$, then a number 
$N(z)$ is attached,  such that $N(z)= \min_X\{-\chi (X)~|~ X\ne S^2 ~\hbox{is compact surface 
representing class} ~z\}$,  where $\chi$ is the Euler number of $X$. The norm is called a {\it Thurston norm}.
The function $N$ is  linear on the $H_2(M)$ and extends to the real homology group 
$H_2(M;{\Bbb R})$ as a pseudo-norm (Thurston \cite{Thu1}).
The Thurston norm is an important homotopy invariant of the manifold $M$, which
can be viewed as a generalization of the genus of a knot.

It is interesting that the abelian groups with a norm
arise in the  context of the  $AF$-algebras (the operator algebras, see Effros \cite{E}). Namely, let 
${\Bbb Z}^k\buildrel\rm A\over\longrightarrow {\Bbb Z}^k
\buildrel\rm A\over\longrightarrow {\Bbb Z}^k\buildrel\rm A\over\longrightarrow \dots$,
be a  stationary dimension group. (We refer the reader to the section 2 for a definition.)  
By $\lambda_A$ one understands the Perron-Frobenius eigenvalue of the 
positive integral matrix $A$.   Let $K={\Bbb Q}(\lambda_A)$
be an algebraic number field of the degree $k$, obtained as 
an extension of the rationals by the algebraic number $\lambda_A$.
The field $K$ is known to be an important invariant of the 
stationary dimension group.  Namely, the triple $(K, \alpha, I)$,
where $\alpha$ is an embedding of the field $K$ and $I$ is the equivalence
class of ideals in the ring of integers of $K$, is a  complete  Morita
invariant of the stationary dimension group (Bratteli, J\o rgensen, 
Kim \& Roush \cite{BJKR}, Effros \cite{E} , Handelman \cite{Han}). 
Denote by $O_K$ the ring of integers of the field $K$ and fix an integral basis 
$\omega_1,\dots,\omega_k$ in $O_K$. Note that $O_K\cong {\Bbb Z}^k$ by the 
Gauss isomorphism.  It is well known that the multiplication by $\alpha\in O_K$ 
induces a linear operator on the  vector space $O_K$: 
$$
\left(
\matrix{
\alpha\omega_1 \cr
\vdots \cr
\alpha\omega_k
}\right)
=
\left(
\matrix{
a_{11} &\dots & a_{1k}\cr
\vdots && \vdots\cr
a_{k1} &\dots & a_{kk}
}\right)
\left(
\matrix{
\omega_1 \cr
\vdots \cr
\omega_k
}\right),
$$
where $a_{ij}$ are the rational integers. Define a function 
${\cal N}: O_K\to {\Bbb Z}$ by the formula 
$\alpha\mapsto a_{11}+\dots+a_{kk}$, where  $\alpha\in O_K$. 
It is not hard to verify,  that ${\cal N}$ is a  linear function, which is independent of 
the choice of the integral basis in $O_K$ (Weyl \cite{W}). Note that the pre-image
${\cal N}^{-1}({\Bbb Z}^+)$ of the semi-group ${\Bbb Z}^+=\{0, 1, 2, \dots\}$ is 
a cone $C\subset {\Bbb Z}^k$.

The aim of this note is to show that the algebraic norm ${\cal N}$ and the Thurston  norm
$N$ are related.  Namely, let $\varphi:X\to X$ be a pseudo-Anosov diffeomorphism
of the compact surface $X$ of the genus $g$ and let ${\cal F}$ 
be the $\varphi$-invariant foliation  on $X$ (Thurston \cite{Thu2}). 
Consider the mapping torus of $\varphi$, i.e. a $3$-dimensional manifold
$M=\left\{ X\times [0,1] ~| \quad (x,0)\mapsto (\varphi(x),1),
\quad x\in X\right\}$. (The reader can recognize $M$ to be the surface bundle over the 
circle with a monodromy $\varphi$.)  Let us  construct the  crossed product $C^*$-algebra 
$A_{\varphi}=C(X)\rtimes_{\varphi}{\Bbb Z}$,
where $C(X)$ is the $C^*$-algebra of the continuous complex-valued 
functions on the surface $X$.  It can be shown, that the $K_0$-group of 
$A_{\varphi}$  is  a stationary dimension group.  Define a map ${\cal N}: H_2(M)\to {\Bbb Z}$
using a natural isomorphism $O_K\cong H_1(X, Sing~{\cal F}; {\Bbb Z})\cong H_2(M)$, where $Sing~{\cal F}$ is the set 
of singular points of the foliation ${\cal F}$.  Our main result is the following theorem.
\begin{thm}\label{thm1}
For every surface bundle $M\to S^1$ with a  pseudo-Anosov monodromy $\varphi$, 
the following is true:  (i) the preimage ${\cal N}^{-1}({\Bbb Z}^+)$ of the semi-group
${\Bbb Z}^+=\{0,1,2,\dots\}$ is a cone $C\subset H_2(M)$ and
(ii) the norm ${\cal N}$ coincides with the Thurston norm $N$ on the cone $C$.
\end{thm}

\medskip\noindent
The structure of the note is the following.  In section 2, the notation
is  introduced. Theorem \ref{thm1} is proved in section 3.

\section{Notation}
This section is a brief introduction to the dimension groups, 
the algebraic number fields and the Thurston norm on the 3-dimensional 
manifolds. We refer the reader to 
O.~Bratteli, P.~E.~T.~ J\o rgensen, K.~H.~ Kim \&  F.~Roush,
(\cite{BJKR}), M.~R\o rdam, F.~Larsen \& 
N.~Laustsen (\cite{RLL}), H.~Weyl (\cite{W}) and 
W.~Thurston (\cite{Thu1}) for a complete account. 
\subsection{The dimension group}
By the $C^*$-algebra one understands a noncommutative Banach
algebra with an involution (\cite{RLL}). Namely, a $C^*$-algebra
$A$ is an algebra over $\Bbb C$ with the norm $a\mapsto ||a||$
and an involution $a\mapsto a^*, a\in A$, such that $A$ is
complete with respect to the norm, and such that 
$||ab||\le ||a||~||b||$ and $||a^*a||=||a||^2$ for every
$a,b\in A$. If $A$ is commutative, then the Gelfand
theorem says that $A$ is isometrically $*$-isomorphic   
to the $C^*$-algebra $C_0(X)$ of the continuous complex-valued
functions on a locally compact Hausdorff space $X$.
For otherwise, $A$ represents a noncommutative topological
space.

\subsubsection{The ordered abelian groups}
Given a $C^*$-algebra, $A$, consider a new $C^*$-algebra $M_n(A)$, 
i.e. the matrix algebra over $A$. There exists a remarkable semi-group, $A^+$,
connected to the set of projections in the algebra
$M_{\infty}=\cup_{n=1}^{\infty} M_n(A)$. Namely, the projections 
$p,q\in M_{\infty}(A)$ are the Murray-von Neumann equivalent $p\sim q$, if they can be 
presented as $p=v^*v$ and $q=vv^*$ for an element $v\in M_{\infty}(A)$.
An equivalence class of the projections is denoted by $[p]$.
The semi-group $A^+$ is defined to be the set of all equivalence
classes of projections in $M_{\infty}(A)$ with the binary operation
$[p]+[q]=[p\oplus q]$. The Grothendieck completion of $A^+$ to an abelian group is called
a {\it $K_0$-group of $A$}.
The functor $A\to K_0(A)$ maps the unital
$C^*$-algebras into a category of the abelian groups, so that
the semi-group $A^+\subset A$ corresponds to a positive
cone  $K_0^+\subset K_0(A)$ and the unit element
$1\in A$ corresponds to an order unit
$[1]\in K_0(A)$. The ordered abelian group $(K_0,K_0^+,[1])$ with the order
unit is called a {\it dimension group} of the $C^*$-algebra $A$.

\subsubsection{The $AF$-algebras}
An $AF$ (approximately finite-dimensional) $C^*$-algebra is defined to
be a  norm closure of an ascending sequence of the finite dimensional
$C^*$-algebras $M_n$'s, where  $M_n$ is a $C^*$-algebra of the $n\times n$ matrices
with the entries in ${\Bbb C}$. Here the index $n=(n_1,\dots,n_k)$ represents
a {\it multi-matrix} $C^*$-algebra $M_n=M_{n_1}\oplus\dots\oplus M_{n_k}$.
Let $M_1\buildrel\rm\varphi_1\over\longrightarrow M_2\buildrel\rm\varphi_2\over\longrightarrow\dots$,
be a chain of the finite-dimensional $C^*$algebras and their homomorphisms. A set-theoretic limit
$A=\lim M_n$ has a natural algebraic structure given by the formula
$a_m+b_k\to a+b$; here $a_m\to a,b_k\to b$ for the
sequences $a_m\in M_m,b_k\in M_k$.
The homeomorphisms of the above (multi-matrix) algebras admit a 
canonical description (Effros \cite{E}).
Suppose that $p,q\in {\Bbb N}$ and $k\in {\Bbb Z}^+$ are such numbers that $kq\le p$. Let us
define a homomorphism $\varphi:M_q\to M_p$ by the formula
$a\longmapsto a\oplus\dots\oplus a\oplus~ 0_h$,
where $p=kq+h$. More generally, if $q=(q_1,\dots,q_s),p=(p_1,\dots,p_r)$
are vectors in ${\Bbb N}^s, {\Bbb N}^r$, respectively, and $\Phi=(\phi_{kl})$ is
a $r\times s$ matrix with the entries in ${\Bbb Z}^+$ such that
$\Phi(q)\le p$, then the homomorphism $\varphi$ is defined by
the formula:
\begin{eqnarray}\label{eq13}
a_1\oplus\dots\oplus a_s &\longrightarrow&
\underbrace{(a_1\oplus a_1\oplus\dots)}_{\phi_{11}}\oplus
\underbrace{(a_2\oplus a_2\oplus\dots)}_{\phi_{12}}\oplus\dots\oplus 0_{h_1}
\nonumber\\
&\oplus&  
\underbrace{(a_1\oplus a_1\oplus\dots)}_{\phi_{21}}\oplus
\underbrace{(a_2\oplus a_2\oplus\dots)}_{\phi_{22}}\oplus\dots\oplus 0_{h_2}
\oplus\dots\nonumber
\end{eqnarray}
where $\Phi(q)+h=p$. We say that $\varphi$  is a {\it
canonical homomorphism} between $M_p$ and $M_q$. 
Any homomorphism $\varphi:M_q\to M_p$ can be rendered canonical
(\cite{E}).

\subsubsection{The Bratteli diagrams}
This graphical presentation of the canonical homomorphism
is called a {\it Bratteli diagram}. Every block
of such a diagram is a bipartite graph with the $r\times s$ matrix
$\Phi=(\phi_{kl})$. In general,
the Bratteli diagram is given by the vertex set $V$ and the edge set $E$, such that
$V$ is an infinite disjoint union $V_1\sqcup V_2\sqcup\dots$, 
where each $V_i$ has
a cardinality $n$. Any pair $V_{i-1},V_i$ defines a non-empty set
$E_i\subset E$ of the edges with a pair of the range and the source functions
$r,s$,  such that $r(E_i)\subseteq V_i$ and $s(E_i)\subseteq V_{i-1}$.
The non-negative integral matrix of the incidences
$(\phi_{ij})$  shows how many edges
are drawn  between the $k$-th vertex in the row $V_{i-1}$ and $l$-th vertex
in the row $V_i$. A Bratteli diagram is called {\it stationary}, if $(\phi_{kl})$ is a
constant  matrix for all $i=1,\dots,\infty$.

\subsection{The number fields}
Let ${\Bbb Q}$ be the field of the rational numbers. Let
$\alpha\not\in {\Bbb Q}$ be an algebraic number over ${\Bbb Q}$, 
i.e. root of the polynomial equation
$a_nx^n+a_{n-1}x^{n-1}+\dots+a_0=0, \qquad a_n\ne 0$,
where $a_i\in {\Bbb Q}$. An algebraic extension of 
the degree $n$ is a minimal field $K=K(\alpha)$, which
contains both ${\Bbb Q}$ and $\alpha$. Note that the coefficients
$a_i$ can be assumed integer. If $K$
is an algebraic extension of the  degree $n$ over ${\Bbb Q}$,
then $K$ is isomorphic to the $n$-dimensional vector space
(over ${\Bbb Q}$) with the basis vectors 
$\{1,\alpha,\dots,\alpha^{n-1}\}$ (Pollard \cite{P}).

\subsubsection{The ring of integers}
Let $K$ be an algebraic extension of the degree $n$ over
${\Bbb Q}$. The element $\tau\in K$ is called 
{\it algebraic integer} if there exits a monic
polynomial $\tau^n+a_{n-1}\tau^{n-1}+\dots+a_0=0$,
where $a_i\in {\Bbb Z}$. It can be easily
verified that the sum and the product of two algebraic integers
is an algebraic integer. The (commutative) ring $O_K\subset K$ 
is called a {\it ring of integers}. The elements of the subring
${\Bbb Z}\subset O_K$ are called the rational integers.  
An {\it integral basis} is a collection $\omega_1,\dots,\omega_n$ of the elements 
of $O_K$,  whose linear span over the rational integers is equal
to $O_K$. An integral basis exists for any finite extension
and therefore $O_K$ is isomorphic to the integral lattice
${\Bbb Z}^n$,  where $n$ is the degree of the field $K$
(Weyl \cite{W}).

\subsubsection{The trace of an algebraic number}
Let $K$ be a number field of the degree $n$ over ${\Bbb Q}$. 
There exists $n$ isomorphic embeddings (monomorphisms)
$K\to {\Bbb C}$ (McCarthy \cite{M}). We denote them
by $\sigma_1,\dots,\sigma_n$. If $\alpha\in K$ then
one defines a {\it trace}  by the formula
\footnote{This notion of ${\cal N}$ is equivalent 
to given in section 1.}
: ${\cal N}(\alpha)=\sigma_1(\alpha)+\dots+\sigma_n(\alpha)$.
When $\alpha$ is an algebraic integer, then ${\cal N}(\alpha)$ is 
a rational integer. If $p,q\in {\Bbb Z}$, then 
${\cal N}(p\alpha+q\beta)=p{\cal N}(\alpha)+q{\cal N}(\beta)$,
for all $\alpha,\beta\in K$. It is not hard to see
that the above formula establishes a homomorphism
${\cal N}:O_K\cong {\Bbb Z}^n\to {\Bbb Z}$, which does not 
depend on the choice of an integral basis in $O_K$
(McCarthy \cite{M}, Weyl \cite{W}).

\subsection{The 3-manifolds}
Let $M$ be a compact oriented $3$-manifold. Suppose that
the second homology group $H_2(M;{\Bbb Z})$ is non-trivial.
There exists a linear mapping of the group into the set
of the positive integers,  which is given by the following construction
of W.~P.~ Thurston (\cite{Thu1}).

\subsubsection{The Thurston norm}
Let $X$ be a connected surface of the genus $g\ge0$. Denote by $\chi(X)$
its Euler characteristic, i.e. an integer number $2-2g$.
The negative part of $\chi(X)$ is defined as 
$\chi_-(X)=\max\{0,-\chi(X)\}$.
If $X$ is not connected, one introduces $\chi_-(X)$
as the sum of the negative parts of the connected components 
of $X$. For a cycle $z\in H_2(M;{\Bbb Z})$,  consider
a non-negative integer
$N(z)=\inf\{\chi_-(X)~|~X ~\hbox{is an embedded surface 
representing}~z\}$. 
The $N(z)$ is called a {\it Thurston norm}.  Given two such cycles
$z_1$ and $z_2$, let $X_1$ and $X_2$ be the surfaces representing
them. There exists a unique way to mend $X_1$
and $X_2$ together to obtain new embedded surface $X$, 
such that $\chi_-(X)=\chi_-(X_1)+\chi_-(X_2)$ (Thurston
\cite{Thu1}). Thus, $N(z)=N(z_1)+N(z_2)$ extends
linearly to the entire group $H_2(M;{\Bbb Z})$.

\subsubsection{The pseudo-Anosov diffeomorphisms}
Let $X$ be a surface of the genus $g\ge 2$. Denote
by Mod~$X=Diff~X/Diff_0X$ the mapping class group
of $X$, i.e. the group of the isotopy classes of the orientation 
preserving diffeomorphisms of $X$. The following classification
of Mod~$X$ is due to J.~Nielsen and W.~P.~Thurston. 
\begin{lem}\label{lm1}
{\bf (\cite{Thu2})}
 Any diffeomorphism $\varphi\in$ Mod $X$    
is isotopic to a diffeomorphism $\varphi'$, such that    
 either  (i) $\varphi'$ has finite order, or    
(ii) $\varphi'$ is pseudo-Anosov (non-periodic) diffeomorphism, or    
(iii) $\varphi'$ is reducible by a system of curves $\Gamma$ surrounded    
 by small tubular neighborhoods $N(\Gamma)$, such that on    
 $M\backslash N(\Gamma)$ $\varphi'$ satisfies either (i) or (ii).    
\end{lem}

\subsubsection{The singularity data}
Let $\varphi$ be a pseudo-Anosov diffeomorphism. There 
exists a pair of $\varphi$-invariant measured foliations 
${\cal F_{\pm}}$ on $X$, such that $\varphi$  expands  along 
${\cal F}_+$  and  contracts along ${\cal F}_-$ with dilatation
factor $\mu>1$ (Thurston \cite{Thu2}).
${\cal F}_+$ and ${\cal F}$ are mutually transversal and have common set
of singular points, which are saddle points with $n\ge3$ prongs. 
For brevity, we let  ${\cal F}={\cal F}_+$. 
Recall that the index of $n$-prong saddle $s_n$ is 
$-{1\over 2}(n-2)$. Therefore 
$\sum_{s_n\in Sing~{\cal F}}{n-2\over 2}= 2g-2$,
where $g$ is the genus of surface $X$. If $m=|Sing~{\cal F}|$
is the total number of the singular points of ${\cal F}$, then
$1\le  m\le 4g-4$,
where the minimum  is attained by
a unique saddle $s_{4g-2}$ and maximum by the set  
 $\{s_3,s_3,\dots,s_3\}$ of $4g-4$ saddles.
We refer to the set $\{s_{i_1},\dots,s_{i_m}\}$ as a {\it singularity
data} of $\cal F$.

\subsubsection{The mapping tori}
Let $\varphi:X\to X$ be a diffeomorphism of the surface
$X$.  One can obtain 3-dimensional manifolds $M=M(\varphi)$
by the formula 
$M=\left\{ X\times [0,1] ~|  ~(x,0)\mapsto (\varphi(x),1),
~x\in X\right\}$.
The manifold $M$ is called a {\it mapping torus}.
It is not hard to see that $M$ is a mapping torus
if and only if $M\to S^1$ is a  fibre bundle over $S^1$ with the monodromy $\varphi$. 
If the diffeomorphism
$\varphi\in Mod~X$ is of a  finite order, then $M$ will be a Seifert
manifold.  In the case when $\varphi$ is pseudo-Anosov,
the following result due to W.~P.~Thurston is true.
\begin{lem}\label{lm2}
{\bf (\cite{Thu3})}
The mapping torus $M$ admits a hyperbolic structure, 
if and only if the diffeomorphism $\varphi$ is pseudo-Anosov.
\end{lem}

\subsubsection{The second homology group of the mapping torus}
Let $\varphi: X\to X$ be a pseudo-Anosov diffeomorphism
of genus $g\ge 2$ surface. Let $Sing~{\cal F}$ be a
finite set of the singularities of the $\varphi$-invariant foliation
$\cal F$. The relative homology group $H_1(X, Sing {\cal F};{\Bbb Z})$
is a torsion-free of the rank
$k=2g +|Sing~{\cal F}|-1$,
where $|Sing~{\cal F}|$ is the cardinality of the set $Sing~{\cal F}$. 
Let $M$ be the mapping torus of $\varphi$. The 2-cycles
of $M$ are generated by the 1-cycles of $X\backslash Sing~{\cal F}$.
Indeed, let $C\in H_1(X, Sing {\cal F};{\Bbb Z})$ and $\varphi(C)$ its image
under the diffeomorphism $\varphi$. Let $X-C$ be a copy of the surface $X$,
which is cut along the closed curve $C$. Similarly, let $X-\varphi(C)$
be a copy of $X$ with a cut along $\varphi(C)$. The surface $X=(X-C)\cup (X-\varphi(C))$,
glued along the action of $\varphi$, belongs to the group $H_2(M)$,
and any non-torsion element  of $H_2(M)$ can be obtained in such a way.   
Therefore, $rank ~H_2(M)=2g +m-1$, where $m=|Sing~{\cal F}|$.

\subsubsection{The Thurston norm of the surface bundle}
It is interesting to relate the Thurston norm $N$
on the second homology with the geometry of
$M$. It turns out,  that in this case $N$ can be expressed
in terms of the Euler classes of the plane bundle tangent
to the fibres of $M$.  
Namely, let $N^*$ be  a dual  Thurston norm defined 
on the first homology $H_1(M;{\Bbb Z})$ by the formula
$N^*(z)=\sup_{N(u)\le 1}u(z)$,
where $u\in Hom (H_1(M;{\Bbb Z}), {\Bbb Z})\simeq H^1(M;{\Bbb Z})$ 
and $N$ is the Thurston norm on $H^1(M;{\Bbb Z})\simeq H_2(M;{\Bbb Z})$
induced by the Poincar\'e duality. Then the following lemma
is true.
\begin{lem}\label{lm3}
{\bf (\cite{Fri}), (\cite{Thu1})}
Let $\tau$ be a subbundle of the tangent
bundle $TM$ consisting of the 2-planes tangent
to the fiber $X$ of the fibration $M\to S^1$. 
Let $e(\tau)\in H^2(M;{\Bbb Z})$ be the Euler
class of $\tau$, i.e. first obstruction 
to the cross-section of bundle $\tau$. Then
(i) the norm $N^*: H_1(M;{\Bbb Z})\to {\Bbb Z}^+$
is induced by the cocycle $e(\tau)$, i.e. 
$N^*(z)=\left|\int_z e(\tau)\right|$;
(ii) the set of the (de Rham) cohomology classes
$H^1(M;{\Bbb R})$, which is  representable by the closed non-singular 
differential 1-forms on $M$ is a maximal cone 
$C\subset H^1(M;{\Bbb R})$, where the Thurston norm $N$ 
can be extended linearly. 
\end{lem}

\section{Proof of theorem 1}
(i) Let ${\cal N}:H_2(M)\to {\Bbb Z}$ be a linear mapping. 
To show that the set $C=\{z\in H_2(M)~|~{\cal N}(z)>0\}$
is a cone, we have to establish that 

\medskip
(1) $z\in C, c>0$ implies $cz\in C$

\smallskip
(2) $z_1,z_2\in C$ implies $z_1+z_2\in C$.

\medskip\noindent
Indeed, since ${\cal N}$ is linear,  ${\cal N}(cz)=c{\cal N}(z)$,
where ${\cal N}(z)>0$ by the assumption. Therefore,
$c{\cal N}(z)>0$ and the item (1) follows. Similarly,
in the item (2), by the linearity of ${\cal N}$,  we have ${\cal N}(z_1+z_2)=
{\cal N}(z_1)+{\cal N}(z_2)>0$, since ${\cal N}(z_1)>0,
{\cal N}(z_2)>0$ by the assumption. The item (i) is proved.

\medskip
(ii) The proof of item (ii) is based on a lemma
of D.~Gabai (\cite{Gab}). Roughly speaking, we shall 
 estimate the Gromov norm (a simplicial norm) of $H_2(M)$,
rather than the Thurston norm itself. This approach gives
a technical advantage, because the group of the 2-chains in
$M$ has a natural abelian structure. Next we use the Gabai
lemma to evaluate the two norms. The Gromov  norm was introduced
and studied in (\cite{Gro}).

Let $M$ be a compact manifold and $z\in H_2(M)$
be an element of the  second homology
group of $M$. A {\it Gromov norm} $g(z)$ is given by the
formula
$g(z)=\inf_{Z\in [z]}\bigl\{\sum |a_i|~:~ Z=\sum a_i\sigma_i
\bigl\}$,
where $[z]$ is the homology class of the 2-chains and
$\sigma_1,\dots,\sigma_n$ is a basis of the
simplicial decomposition of $M$. The following lemma is true.
\begin{lem}\label{lm4}
Suppose $M$ is a compact oriented 3-manifold.
Then the Thurston norm $N$ and Gromov norm $g$
are related by the formula $N(z)={1\over 2}g(z)$,
for each $z\in H_2(M)$ in the domain of definition
of the two norms.  
\end{lem}
{\it Proof.} The proof can be found
in (\cite{Gab}), Corollary 6.18. For the sake of clarity, 
let us outline  the main idea. First, notice
that if $M$ is a hyperbolic $k$-manifold and $[M]$ its homology class, 
then we have Gromov's formula $g([M])=Vol~M~/~Vol~\sigma$, 
where $\sigma$ is the largest hyperbolic $k$-simplex 
(Gromov \cite{Gro}).  Thus, for the connected surface $X$,
one has $g([X])={2\pi|\chi(X)|\over\pi}=2N([X])$.
The formula extends to the case $X$ with more than one
connected component and requires the singular norms
$x_s$ in this case, see Gabai (\cite{Gab}). Eventually,
it can be shown,  that $x_s=N$ and one gets inequality
$g\le 2N$.

To prove the inequality  $g\ge 2N$, let $z\in H_2(M)$ 
and $Z\in [z]$ be a 2-cycle $Z=\sum a_i\sigma_i$, where 
$a_i\in {\Bbb Z}$. By pasting the singular simplices, one can obtain a 
proper map $f:X\to M$,  such that $[f(X)]=z$. 
By Gromov's formula for the hyperbolic volumes
$2N([X])\le \sum|a_i|= g([X])$.
Lemma \ref{lm4} follows.
$\square$

\medskip\noindent
Fix a simplicial basis $\sigma_1,\dots,\sigma_n$ 
in the group $\cal C$ of the 2-chains of the regular 
triangulation of $M$. Consider a subset of $\cal C$,
given by the formula
$K_{\cal C}=\{Z\in {\cal C}~|~Z=\sum_{i=1}^na_i\sigma_i,
\quad a_i>0\}$.
Denote by $\widetilde {\cal N}: H_2(M)\to {\Bbb Z}^+$
the mapping $z\mapsto g(z)$ given by the Gromov norm.
It is not hard to verify,  that $\widetilde {\cal N}$ is
linear on the $K_{\cal C}$ and the $K_{\cal C}$ is a
cone in $\cal C$. Consider the following commutative
diagram of the linear mappings:

\bigskip
\begin{picture}(400,90)(50,10)
\put(130,80){$K_{\cal C}$}
\put(150,83){\vector(3,0){60}}
\put(220,80){$C$}
\put(175,90){$\Sigma$}
\put(135,70){\vector(0,-3){35}}
\put(225,70){\vector(0,-3){35}}
\put(130,20){${\Bbb Z}^+$}
\put(150,23){\vector(3,0){60}}
\put(220,20){${\Bbb Z}^+$}
\put(175,30){$w$}
\put(120,50){$\widetilde {\cal N}$}
\put(232,50){${\cal N}$}
\end{picture}

\medskip\noindent
(Here $w$ is a  doubling map,  acting by the formula
$z\mapsto 2z$). The  map $\Sigma$ on
the diagram is linear. The $\Sigma$ establishes
a bijection between the bases $\{\sigma_i\}$
and $\{\omega_i\}$: $\Sigma(\sigma_i)=\omega_i$.  
It remains to apply lemma \ref{lm4}. The item (ii) of
theorem \ref{thm1} follows.
$\square$

\medskip\noindent
{\sf Acknowledgments.}  I wish to thank O.~Bratteli, 
D.~Calegari, M.~Felisatti, Y.~Minsky and V.~Turaev for 
helpful comments.



\vskip1cm

\textsc{The Fields Institute for Mathematical Sciences, Toronto, ON, Canada,  
E-mail:} {\sf igor.v.nikolaev@gmail.com} 

\smallskip
{\it Current address: 101-315 Holmwood Ave., Ottawa, ON, Canada, K1S 2R2}


\begin{thebibliography}{100}
\bibitem{BJKR}
O.~Bratteli, P.~E.~T.~ J\o rgensen, K.~H.~ Kim, and F.~Roush,  
Computation of isomorphism invariants for stationary dimension groups. Ergodic Theory Dynam. Systems 22 (2002), 
 99-127.


\bibitem{E}
E.~G.~Effros, Dimensions and $C^*$-Algebras, in: Conf. Board of the Math.
Sciences No.46, AMS (1981).



\bibitem{Fri} 
 D.~Fried, Fibrations over $S^1$ with pseudo-Anosov monodromy,
 in: A.~Fathi, F.~Laudenbach and V.~Po\'enaru,
 Travaux de Thurston sur les Surfaces, Expos\'e 14, Ast\'erisque 
  66-67 (1979), 251-266.


\bibitem{Gab}
D.~Gabai, Foliations and the topology of 3-manifolds,
J. Differential Geometry 18 (1983), 445-503. 



\bibitem{Gro}
M.~Gromov, Volume and bounded cohomology, Inst. Hautes \'Etudes Sci. 
Publ. Math. 56 (1982), 5-99.



\bibitem{Han}
D.~Handelman, Positive matrices and dimension groups affiliated
to $C^*$-algebras and topological Markov chains, J. Operator
Theory 6 (1981), 55-74. 


\bibitem{M}
P.~J.~McCarthy, Algebraic Extensions of Fields,
Blaisdell Publishing Company, 1966.



\bibitem{P}
H.~Pollard, The Theory of Algebraic Numbers, Carus mathematical
monographs, 9. 


\bibitem{RLL}
M.~R\o rdam, F.~Larsen and N.~Laustsen, An introduction to $K$-theory 
for $C^*$-algebras. London Mathematical Society Student Texts, 49. 
Cambridge University Press, Cambridge, 2000. xii+242 pp. ISBN: 0-521-78334-8; 
0-521-78944-3



\bibitem{Thu1}
W.~P.~Thurston, A norm for the homology of 3-manifolds,
Mem. Amer. Math. Soc. 339 (1986), 99-130.


\bibitem{Thu2}
W.~P.~Thurston, On the geometry and dynamics of diffeomorphisms 
of surfaces, Bull. Amer Math. Soc. 19 (1988), 417-431.


\bibitem{Thu3}
W.~P.~Thurston, Hyperbolic structures on 3-manifolds, II:
surface groups and 3-manifolds which fiber over the circle,
arXiv:math.GT/9801045 v1.



\bibitem{W}
H.~Weyl, Algebraic Theory of Numbers, Annals of Math. Studies 1, 
Princeton Univ. Press, 1940. 


\end{thebibliography}
\end{document}